\documentclass[11pt]{amsart}

\usepackage{amscd}
\usepackage{amsmath, amssymb}
\usepackage{amsfonts}
\usepackage{color}
\newcommand{\de}{\partial}

\newcommand{\ddbar}{\sqrt{-1} \partial \overline{\partial}}

\newcommand{\ov}[1]{\overline{#1}}

\newcommand{\ti}[1]{\tilde{#1}}
\newcommand{\vp}{\varphi}

\newcommand{\ve}{\varepsilon}

\renewcommand{\bar}{\overline}

\renewcommand{\leq}{\leqslant}
\renewcommand{\geq}{\geqslant}

\numberwithin{equation}{section}

\begin{document}
\newtheorem{claim}{Claim}
\newtheorem{theorem}{Theorem}[section]
\newtheorem{lemma}[theorem]{Lemma}
\newtheorem{corollary}[theorem]{Corollary}
\newtheorem{proposition}[theorem]{Proposition}
\newtheorem{conjecture}[theorem]{Conjecture}

\theoremstyle{definition}
\newtheorem{remark}[theorem]{Remark}
\newtheorem{question}[theorem]{Question}

\title{A singular Demailly-P\u{a}un theorem}
\author[T.C. Collins]{Tristan C. Collins}
\address{Department of Mathematics, Harvard University, 1 Oxford Street, Cambridge, MA 02138}
\email{tcollins@math.harvard.edu}
\author[V. Tosatti]{Valentino Tosatti}
\address{Department of Mathematics, Northwestern University, 2033 Sheridan Road, Evanston, IL 60208}
\email{tosatti@math.northwestern.edu}
\begin{abstract}
We give a numerical characterization of the K\"ahler cone of a possibly singular compact analytic variety which is embedded in a smooth ambient space.
\end{abstract}

\maketitle

\section{Introduction}
The classical Nakai-Moishezon ampleness criterion (see e.g. \cite{Kl} and references therein) characterizes ample line bundles on a projective variety as those which have positive intersection against all subvarieties. This was later extended to $\mathbb{R}$-divisors by Campana-Peternell \cite{CP}. In a groundbreaking paper, Demailly and P\u{a}un \cite{DP} proved a vast generalization of this result, which holds for all real $(1,1)$ classes on a compact K\"ahler manifold. More precisely, they proved that the K\"ahler cone of a compact K\"ahler manifold is one of the connected components of the positive cone, consisting of classes which have positive intersection against all analytic subvarieties. Very recently, a new proof of this theorem was obtained by combining the main result of our previous work \cite{CT} with a result of Chiose \cite{Ch}.

In this note, we prove an extension of the Demailly-P\u{a}un theorem \cite{DP} to singular varieties which are embedded in a smooth ambient space. A $(1,1)$ class on the variety is just taken to be the restriction of a $(1,1)$ class from the ambient space, and such a class is K\"ahler if it is so in a neighborhood of the variety inside the ambient space. This is in fact equivalent to the more intrinsic definition of a K\"ahler class on a compact analytic space as given for example in \cite{Va}, as shown by P\u{a}un \cite{Pa}, and this allows us to avoid discussing these more technical notions. 
With these observations in mind, our main theorem is the following:

\begin{theorem}\label{main}
Let $(M,\omega)$ be a smooth (but possibly noncompact and incomplete) K\"ahler manifold, and $E\subset M$ be a compact analytic subvariety. Let $\alpha$ be a closed smooth real $(1,1)$ form on $M$ such that
\begin{equation}\label{pos}
\int_V \alpha^{k}\wedge\omega^{\dim V-k}>0,
\end{equation}
for all positive-dimensional irreducible analytic subvarieties $V\subset E$, and for all $1\leq k\leq\dim V$. Then there exist an open neighborhood $U$ of $E$ in $M$ and a smooth function $\vp:U\to\mathbb{R}$ such that $\alpha+\ddbar \vp$ is a K\"ahler metric on $U$.
If $M$ is an open subset of the regular locus of some projective variety, then the inequalities
\begin{equation}\label{pos1}
\int_V \alpha^{\dim V}>0,
\end{equation}
for all $V$ as above suffice to reach the same conclusion.
\end{theorem}

This theorem answers a question that was posed to us by R.J. Conlon and H.-J. Hein, in relation to their paper \cite{CH} (see also \cite[1.3.5]{CH1}). Applications of this result to the study of the K\"ahler cone of asymptotically conical Calabi-Yau manifolds will appear in a forthcoming revision of \cite{CH}.

The main tools we use are the Demailly-P\u{a}un theorem itself, for smooth compact K\"ahler manifolds, and our recent theorem \cite{CT} which shows that the non-K\"ahler locus of a nef and big class on a compact complex manifold equals the null locus of the class. The idea is to work by induction on the dimension on $E$ (as in \cite{DP}), and to prove the result by working on a resolution of singularities (as in \cite{CT}). This way we avoid any technical discussion of currents on singular analytic spaces. 

In future work, we hope to address the extension of the 
Demailly-P\u{a}un theorem \cite{DP} as well as the main result of our previous work \cite{CT} to general compact K\"ahler (reduced and irreducible) analytic spaces.\\

{\bf Acknowledgments. }The first-named author is supported in part by NSF grant DMS-1506652, and the second-named author by NSF grant DMS-1308988 and by a Sloan Research Fellowship. Part of this work was carried out while the second-named author was visiting the Yau Mathematical Sciences Center of Tsinghua University in Beijing, which he would like to thank for the hospitality. We are grateful to R.J. Conlon and H.-J. Hein for asking the question which is answered in this paper, and to H.-J. Hein for some very helpful discussions.

\section{Proof of Theorem \ref{main}}
This section contains the proof of Theorem \ref{main}.

Clearly we may assume that no component of $E$ is zero-dimensional, since for those the result is trivial.

Let us first assume that $E$ is irreducible and $1$-dimensional. Let $\nu:\ti{M}\to M$ be an embedded resolution of singularities of $E\subset M$, so that $\ti{M}$ is smooth, connected and K\"ahler, and the proper transform $\ti{E}$ of $E$ is a smooth compact Riemann surface. We will also write $\nu:\ti{E}\to E$ for the induced map, so that $\nu^*\alpha$ is a smooth closed real $(1,1)$ form with $\int_{\ti{E}}\nu^*\alpha>0$. Therefore the class $[\nu^*\alpha]$ on $\ti{E}$ is K\"ahler, and we can find a smooth function $\psi$ on $\tilde{E}$ such that $\nu^*\alpha+\ddbar \psi>0$ on $\ti{E}$. It is elementary to find an open neighborhood $\ti{U}$ of $\ti{E}$ in $\ti{M}$ and a smooth extension of $\psi$ to $\ti{U}$ (still denoted by $\psi$) such that $\nu^*\alpha+\ddbar\psi>0$ on $\ti{U}$ (see e.g. \cite[Lemme 1, p.416]{Pa}). Note that $U=\nu(\ti{U})\backslash E_{sing}$ is an open neighborhood of $E_{reg}$ inside $M$, but in general it is not the case that $\nu(\ti{U})$ is an open neighborhood of $E$ inside $M$, because it may ``pinch off'' near $E_{sing}$. Furthermore, even if $\nu(\ti{U})$ happens to be an open neighborhood of $E$, the pushforward function $\nu_*\psi$ is not well-defined wherever different branches of $\ti{E}$ come together under the map $\nu$. Therefore, we need to work a bit harder to achieve our goal.

Let $\{p_1,\dots,p_N\}\subset \tilde{E}$ be the exceptional locus of $\nu$ intersected with $\ti{E}$, so that $\{\nu(p_1),\dots,\nu(p_N)\}$ equals the singular set of $E$. For each point $p_j$ we add to $\psi$ a function of the form $\ve\theta(z) \log|z-p_j|$, where $\ve>0$ is small enough, where $z=(z_1,\dots,z_{N})$ are local coordinates for $\ti{M}$ near $p_j$, and $\theta$ is a smooth cutoff function supported in a small neighborhood of $p_j$ in $\ti{M}$, so that
we obtain a new function $\ti{\psi}$, which is smooth away from the $p_j$'s and goes to $-\infty$ there, and such that $\nu^*\alpha+\ddbar\ti{\psi}$ is a K\"ahler current on $\tilde{U}$.

Then the smooth function $\hat{\psi}=\nu_*\ti{\psi}$ on $U$ satisfies $\alpha+\ddbar\hat{\psi}>0$, but we are not done yet because $U$ does not contain the singular points of $E$. Let $\{\nu(p_1),\dots,\nu(p_k)\}$ be all the singular points of $E$ (so $k\leq N$), and fix charts $U_j$ for $M$ centered at $\nu(p_j)$ for $1\leq j\leq k$, with coordinates so that each $U_j$ is the Euclidean ball of radius $2$. Call $U'_j$ the Euclidean ball of radius $1$ in these coordinates, and let $A$ be the minimum of $\hat{\psi}$ on the compact set
$$\bigcup_{j=1}^k \ov{(\de U'_j)\cap U},$$
which is a finite number because $\hat{\psi}$ is smooth there. Choose a large constant $B>0$ such that on each $U_j$ we have $\alpha+B\ddbar |z|^2>0$.
On $U\cap U_j$ then we have that $\hat{\psi}$ and $B|z|^2+A-B-1$ are both strictly $\alpha$-plurisubharmonic, with $\hat{\psi}$ approaching $-\infty$ at the center of the ball $U_j$, and with $\hat{\psi}>B|z|^2+A-B-1$ on a neighborhood of $(\de U'_j)\cap U$. If $\widetilde{\max}$ denotes a regularized maximum function  (see, e.g. \cite[I.5.18]{Demb}), then
$$\psi_g=\widetilde{\max}(\hat{\psi}, B|z|^2+A-B-1)$$
is smooth and strictly $\alpha$-plurisubharmonic on $U_j\cap U$, it equals $\hat{\psi}$ in a neighborhood of $(\de U'_j)\cap U$, and it equals $B|z|^2+A-B-1$ as we approach the origin. Therefore the function $\psi_g$ trivially glues to $\hat{\psi}$ outside $U'_j$, and we can extend it to be equal to  $B|z|^2+A-B-1$ in a small neighborhood of the origin in $U_j$. Repeating this construction for all $j$, and gluing each of them to $\hat{\psi}$, we finally obtain an open neighborhood $\bar{U}$ of $E$ in $M$ and a smooth function $\vp$ on $\bar{U}$ such that
$\alpha+\ddbar\vp$ is a K\"ahler metric on $\bar{U}$, as required.

Next, we assume that $E$ has pure dimension $1$, but need not be irreducible anymore. Then, writing $E=\cup_j E_j$ with $E_j$ irreducible, we can apply the result to each $E_j$ and obtain $U_j,\vp_j$ as above, and ``glue'' them all together using \cite[Lemme, p.419]{Pa}, and obtain the desired K\"ahler potential $\vp$ on some neighborhood $\bar{U}$ of $E$.

We now deal with the general case, by induction on $\dim E$ (which is by definition the max of the dimensions of the irreducible components of $E$). The base of the induction is what we have just proved. For the induction step, let $\dim E=n$ and assume the result holds in all dimensions $<n$. As we just did, it is enough to prove the theorem in the case when $E$ is irreducible, since if there are several components then we work on each one separately, and in the end glue the resulting metrics as before. So we will assume that $E$ is irreducible.
Take $\nu:\ti{M}\to M$ to be an embedded resolution of singularities of $E\subset M$, obtained as a composition of blowups with smooth centers, so that $\ti{M}$ is smooth and K\"ahler, and the proper transform $\ti{E}$ of $E$ is smooth.

Then $\nu^*\alpha$ is a smooth closed real $(1,1)$ form on $\ti{E}$, and we claim that its class $[\nu^*\alpha]$ on $\ti{E}$ is nef. If assume that $M$ is an open subset of the regular locus of some projective variety, then this holds because we have $\int_{V}(\nu^*\alpha)^{\dim V}\geq 0$ for all positive-dimensional irreducible subvarieties $V$ in $\ti{E}$ (using \eqref{pos}), and so \cite[Theorem 4.5(ii)]{DP} gives that the class $[\nu^*\alpha]$ on $\ti{E}$ is nef. However, in our general setup (where there may be no projective compactification), to use \cite[Theorem 4.3(ii)]{DP} we would have to check instead that
$$\int_{V}\nu^*\alpha^{k}\wedge\tilde{\omega}^{\dim V-k}\geq 0,$$
for all positive-dimensional irreducible subvarieties $V\subset\ti{E}$, for some K\"ahler form $\ti{\omega}$ on $\ti{E}$ and for all $1\leq k\leq \dim V$, and it does not seem easy to check this directly. Instead, we argue as follows. We have
$$\int_{\ti{E}}\nu^*(\alpha^k \wedge \omega^{\dim E-k})>0,$$
for $1\leq k\leq \dim E$, because $\nu:\ti{E}\to E$ is a modification, and using \eqref{pos}. Since the class $[\nu^*\omega]$ is nef on $\ti{E}$, we can find K\"ahler classes on $\ti{E}$ arbitrarily close to it, and therefore there exists a K\"ahler metric $\ti{\omega}$ on $\ti{E}$ such that
\begin{equation}\label{todo2}
\int_{\ti{E}}\nu^*\alpha^k \wedge \ti{\omega}^{\dim E-k}>0,
\end{equation}
for $1\leq k\leq \dim E$.
Now for $t\geq 0$ sufficiently large, the class $[\nu^*\alpha+t\ti{\omega}]$ is K\"ahler on $\ti{E}$. Let $t_0$ be the minimum value of $t$ such that the class
$[\nu^*\alpha+t\ti{\omega}]$ is nef on $\ti{E}$, and suppose for a contradiction that $t_0>0$. By definition the class $[\nu^*\alpha+t_0\ti{\omega}]$ is not K\"ahler on $\ti{E}$.
Thanks to \cite[Theorem 0.1]{DP}, there exists a positive-dimensional irreducible analytic subvariety $V\subset \ti{E}$, such that
\begin{equation}\label{todo3}
\int_{V}(\nu^*\alpha+t_0\ti{\omega})^{\dim V}=0,
\end{equation}
since if we had strict positivity for all such $V$ then the class $[\nu^*\alpha+t_0\ti{\omega}]$  would be K\"ahler.
Also $V$ must be properly contained in $\ti{E}$, because we have
$$\int_{\ti{E}}(\nu^*\alpha+t_0\ti{\omega})^{\dim E}>0,$$
by \eqref{todo2}.
Then $\nu(V)$ is an irreducible analytic subvariety of $E$ (possibly a point), of dimension strictly less than $\dim E$, and with the same positivity property \eqref{pos}, so by induction we can find an open neighborhood $W$ of $\nu(V)$ in $M$ and a smooth function $\eta$ on $W$ such that $\alpha+\ddbar\eta>0$. Therefore, in the open neighborhood $\nu^{-1}(W)$ of $V$ the smooth function $\nu^*\eta$ satisfies $\nu^*\alpha+\ddbar(\nu^*\eta)\geq 0$. Since $\ti{\omega}$ is K\"ahler on $\ti{E}$ and $t_0>0$, this implies that
$$\int_{V}(\nu^*\alpha+t_0\ti{\omega})^{\dim V}>0,$$
contradicting \eqref{todo3}. Therefore we must have $t_0\leq 0$, and so the class $[\nu^*\alpha]$ is indeed nef on $\ti{E}$.

This proves our claim that the class $[\nu^*\alpha]$ is nef on $\ti{E}$,
and since $$\int_{\ti{E}}(\nu^*\alpha)^{\dim E}=\int_E\alpha^{\dim E}>0,$$ by \eqref{pos}, we can apply \cite[Theorem 2.12]{DP} and see that this class is also big, i.e. it contains a K\"ahler current $\nu^*\alpha+\ddbar\psi$, which we may assume has analytic singularities thanks to Demailly's regularization theorem (see \cite[Theorem 3.2]{DP}). Also, if $V\not\subset  \mathrm{Exc}(\nu)\cap \ti{E}$ then $\nu(V)$ is an irreducible subvariety of $E$ of the same dimension as $V$, and $\nu:V\to\nu(V)$ is bimeromorphic and so we have $\int_{V}(\nu^*\alpha)^{\dim V}=\int_{\nu(V)}\alpha^{\dim V}>0$, thanks to assumption \eqref{pos}. This means that
the null locus of the class $[\nu^*\alpha]$ on $\ti{E}$ is contained in $\mathrm{Exc}(\nu)$, and so using \cite[Theorem 1.1]{CT}, we may choose $\psi$ to be smooth on $\ti{E}\backslash \mathrm{Exc}(\nu)$. We use \cite[Lemma 2.1]{DP} to obtain a quasi-plurisubharmonic function with nontrivial analytic singularities along $\mathrm{Exc}(\nu)$, and add a small multiple of it to $\psi$, to obtain a function $\ti{\psi}$ which is smooth on $\ti{E}\backslash \mathrm{Exc}(\nu)$ and goes to $-\infty$ along $\mathrm{Exc}(\nu)$, and such that $\nu^*\alpha+\ddbar\ti{\psi}$ is a K\"ahler current on $\tilde{E}$ with analytic singularities along $\mathrm{Exc}(\nu)$.

As in the first part of the proof of \cite[Theorem 3.2]{CT}, up to modifying $\ti{\psi}$ slightly (maintaining its same properties) we can find an extension $\ti{\psi}'$ to an open neighborhood $\ti{U}$ of $\ti{E}\backslash \mathrm{Exc}(\nu)$ in $\ti{M}$.

Here are some details for this construction (see \cite{CT} for full details). By a resolution of singularities argument, we construct a modification $\mu:\hat{M}\to\ti{M}$, which is a composition of blowups with smooth centers, such that $\mu(\mathrm{Exc}(\mu))$ is equal to $\mathrm{Exc}(\nu)$, such that the proper transform $\hat{E}$ of $\ti{E}$ is smooth, and the pullback under
$\mu$ of the ideal sheaf on $\ti{E}$ which defines the singularities of the K\"ahler current $\nu^*\alpha+\ddbar\ti{\psi}$ is a principal ideal, supported along a simple normal crossings divisor, which is the restriction to $\hat{E}$ of a simple normal crossings divisor on $\hat{M}$ (which is equal to $\mathrm{Exc}(\mu)$), which has normal crossings with $\hat{E}$. We then cover $\hat{E}$ by finitely many coordinate charts $\{W_j\}$ for $\hat{M}$. To the pullback $\mu^*\ti{\psi}$ we add a small multiple of $\ddbar\log|s|^2_h$, where $s$ defines $\mathrm{Exc}(\mu)$ (and $h$ is chosen suitably), to obtain a strictly $\mu^*\nu^*\alpha$-plurisubharmonic function $\Psi$ on $\hat{E}$, with analytic singularities as before (in particular, smooth away from $\mathrm{Exc}(\mu)$). For each $j$, we then extend $\Psi|_{W_j\cap \hat{E}}$ to a function $\psi_j$ on $W_j$ in an elementary fashion, still preserving strict  $\mu^*\nu^*\alpha$-plurisubharmonicity. Then we use a gluing procedure inspired by a classical method of Richberg \cite{Ri} (see e.g. \cite[Lemma 3.3]{Sm}), but with the extra difficulty that now the functions $\psi_j$ have poles. Nevertheless, arguing exactly as in \cite[Proof of Theorem 3.2]{CT}, we can obtain an open neighborhood $U_1$ of $\hat{E}$ in $\hat{M}$ and a strictly $\mu^*\nu^*\alpha$-plurisubharmonic function $\ti{\Psi}$ on $U_1$, which restricts to $\Psi$ on $\hat{E}$, and is smooth on $\hat{E}\backslash \mathrm{Exc}(\mu)$.

Here we highlight that since $\hat{E}$ is a complex submanifold of a complex manifold, constructing this extension $\ti{\Psi}$ on an open neighborhood $U_1$ of $\hat{E}$ would be standard by Richberg \cite{Ri} if $\Psi$ was smooth (or even just continuous) on $\hat{E}$. On the other hand, if the singularities of $\Psi$ were completely arbitrary, then such an extension would not be possible in general. The key property that saves us here is that the singular locus of $\Psi$ is the intersection with $\hat{E}$ of a simple normal crossings divisor, ${\rm Exc}(\mu)$, in the ambient space $\hat{M}$.

Then we take $\ti{U}=\mu(U_1)$, and $\ti{\psi}'=\mu_*\ti{\Psi}$, which are as required. In particular, $\ti{U}$ is an open neighborhood of $\ti{E}\backslash \mathrm{Exc}(\nu)$ in $\ti{M}$, and $\ti{\psi}'$ is strictly $\nu^*\alpha$-plurisubharmonic, and it is smooth on $\ti{U}\backslash \mathrm{Exc}(\nu)$.

On the open set $U=\nu(\ti{U})\backslash E_{sing}$ (which is a neighborhood of $E_{reg}$ in $M$) we have the smooth function $\hat{\psi}=\nu_*\ti{\psi}'$ with $\alpha+\ddbar (\nu_*\ti{\psi}')$ a smooth K\"ahler metric there, and with $\nu_*\ti{\psi}'$ approaching $-\infty$ along $E_{sing}$. Now $E_{sing}$ is a subvariety of $M$ of dimension strictly less than $n$, with the same positivity property \eqref{pos}, so by induction we can find an open neighborhood $W$ of $E_{sing}$ in $M$ and a smooth function $\hat{\vp}$ on $W$ with $\alpha+\ddbar\hat{\vp}>0$ on $W$. We may also assume that $\hat{\vp}$ is defined on a slightly larger open set, so that it is smooth up to $\de W$.

If we let $A$ be the minimum of $\hat{\psi}$ on the compact set $\overline{(\de W)\cap U}$ and $B$ be the maximum of $\hat{\vp}$ on the same set, then $\hat{\psi}>\hat{\vp}+A-B-1$ holds on a neighborhood of $(\de W)\cap U$. Then
$$\psi_g=\widetilde{\max}(\hat{\psi}, \hat{\vp}+A-B-1)$$
is smooth and strictly $\alpha$-plurisubharmonic on $U\cap W$, equal to $\hat{\psi}$ near $(\de W)\cap U$, and equal to $\hat{\vp}+A-B-1$ as we approach $E_{sing}$. Therefore $\psi_g$ trivially glues to $\hat{\psi}$ outside $W$, and we can extend it to be equal to $\hat{\vp}+A-B-1$ in a neighborhood of $E_{sing}$. In this way we obtain an open neighborhood $\bar{U}$
of $E$ in $M$ and a smooth function $\vp$ on $\bar{U}$ such that
$\alpha+\ddbar\vp$ is a K\"ahler metric on $\bar{U}$, as required.

Lastly, the statement in the projective case follows from the K\"ahler one exactly as in \cite{DP}, by choosing $\omega$ to be the curvature form of a very ample line bundle $L$ on the projective variety which contains $M$ as an open subset, and observing that
$$\int_V \alpha^{k}\wedge\omega^{\dim V-k}=\int_{V\cap H_1\cap \dots\cap H_{\dim V-k}}\alpha^k,$$
for generic members $H_1,\dots,H_{\dim V-k}$ of the linear system $|L|$, so that $V\cap H_1\cap \dots\cap H_{\dim V-k}$ is an irreducible subvariety of dimension $k$.


\begin{thebibliography}{99}
\bibitem{CP} F. Campana, T. Peternell {\em Algebraicity of the ample cone of projective varieties}, J. Reine Angew. Math. {\bf 407} (1990), 160--166. 
\bibitem{Ch} I. Chiose  {\em The K\"ahler rank of compact complex manifolds}, preprint, arXiv:1308.2043.
\bibitem{CT} T. Collins, V. Tosatti {\em K\"ahler currents and null loci}, to appear in Invent. Math.
\bibitem{CH1} R.J. Conlon, H.-J. Hein {\em Asymptotically conical Calabi-Yau manifolds, I}, Duke Math. J. {\bf 162} (2013), no. 15, 2855--2902.
\bibitem{CH} R.J. Conlon, H.-J. Hein {\em Asymptotically conical Calabi-Yau manifolds, III}, preprint, arXiv:1405.7140.
\bibitem{Demb} J.-P. Demailly {\em Complex Analytic and Differential Geometry}, available on the author's webpage.
\bibitem{DP} J.-P. Demailly, M. P\u{a}un \emph{Numerical characterization of the K\"ahler cone of a compact K\"ahler manifold}, Ann. of Math., {\bf 159} (2004), no. 3, 1247--1274.
\bibitem{Kl} S.L. Kleiman {\em  Toward a numerical theory of ampleness}, Ann. of Math. (2) {\bf 84} (1966), 293--344.
\bibitem{Pa} M. P\u{a}un {\em Sur l'effectivit\'e num\'erique des images inverses de fibr\'es en droites}, Math. Ann. {\bf 310} (1998), no. 3, 411--421.
\bibitem{PS} D.H. Phong, J. Sturm {\em On the singularities of the pluricomplex Green's function}, in {\em Advances in analysis. The legacy of Elias M. Stein}, 419--435, Princeton University Press, 2014.
\bibitem{Ri} R. Richberg {\em Stetige streng pseudokonvexe Funktionen}, Math. Ann. {\bf 175} (1968), 257--286.
\bibitem{Sm} P.A.N. Smith {\em Smoothing plurisubharmonic functions on complex spaces}, Math. Ann. {\bf 273} (1986), no. 3, 397--413.
\bibitem{Va} J. Varouchas {\em K\"ahler spaces and proper open morphisms}, Math. Ann. {\bf 283} (1989), no. 1, 13--52.
\end{thebibliography}
\end{document}